\documentclass[11pt]{amsart}
\usepackage{amsmath}
\usepackage{latexsym,ulem}
\usepackage{amssymb}
\usepackage{graphics}

\newtheorem{theorem}{Theorem}[section]

\newtheorem{lemma}{Lemma}[section]

\newtheorem{Definition}{Definition}[section]
\def\qed{{\hfill{\vrule height4pt width3pt depth2pt}}}

   \long\def\comment#1{}

\def\ad#1{\begin{aligned}#1\end{aligned}}  
\def\a#1{\begin{align*}#1\end{align*}} \def\an#1{\begin{align}#1\end{align}}
\def\e#1{\begin{equation}#1\end{equation}} \def\d{\operatorname{div}}
\def\p#1{\begin{pmatrix}#1\end{pmatrix}} 
 
\DeclareMathOperator{\diag}{diag}
\DeclareMathOperator{\sspan}{span}
 \numberwithin{equation}{section}
\numberwithin{table}{section}
\numberwithin{figure}{section}

\def\boxit#1{\vbox{\hrule height1pt \hbox{\vrule width1pt\kern1pt
     #1\kern1pt\vrule width1pt}\hrule height1pt }}
 \def\lab#1{\boxit{\small #1}\label{#1}}
  
 \def\meqref#1{\boxit{\small #1}\eqref{#1}}

   \def\lab#1{\label{#1}}  \def\meqref#1{\eqref{#1}}

\newcommand{\disp}{\displaystyle}

\newcommand{\bq}{\begin{equation}}
\newcommand{\eq}{\end{equation}}

\begin{document}

\title[Superconvergence of simple conforming MFEMs for linear elasticity]
{Superconvergence of simple conforming mixed finite
   elements for linear elasticity on rectangular grids
  in any space dimension
    }
\date{}
\author {Jun Hu}
\address{LMAM and School of Mathematical Sciences, Peking University,
  Beijing 100871, P. R. China.  hujun@math.pku.edu.cn}

\author {Shangyou Zhang}
\address{Department of Mathematical Sciences, University of Delaware,
    Newark, DE 19716, USA.  szhang@udel.edu }

\thanks{The first author was supported by  the NSFC Project 11271035, and  in part by the NSFC Key Project 11031006.}

\begin{abstract}
This paper is to prove superconvergence of   a  family of
  simple conforming mixed finite elements of first order
 for the  linear elasticity problem with the Hellinger--Reissner  variational formulation.
 The analysis is based on three main ingredients: a new interpolation operator, a new expansion method,
 and a new iterative argument for superconvergence analysis.

  \vskip 15pt


 \vskip 15pt


\end{abstract}
\maketitle

\section{Introduction}

This paper investigates superconvergence of  simple conforming
    mixed finite elements \cite{HuManZhang2013}
    for linear  elasticity within the Hellinger-Reissner variational principle.
 It is well--known that it is a challenge problem for stable
    discretizations for this problem, which  results from
    a strong coupling of  the symmetry  requirement  on the  discrete stress tensor
    and the usual stable  conditions for mixed finite element methods.
A lot of efforts, see, for instance,  \cite{Amara-Thomas,
   Arnold-Brezzi-Douglas, Arnold-Douglas-Gupta,  Johnson-Mercier, Morley,
   Stenberg-1, Stenberg-2, Stenberg-3},
    have been devoted to developing stable methods of this problem.
But no stable mixed finite element was found
      in the first four decades \cite{Arnold-Winther-conforming}.
Not until the year 2002, were there some advances
   in this direction.
In \cite{Arnold-Winther-conforming} and
   \cite{Arnold-Awanou-Winther}, a sufficient condition is proposed,
   which states  that a discrete exact sequence
   guarantees the stability of the mixed method.
From then,   conforming mixed finite elements on the
  simplical and rectangular triangulations have been  constructed
    \cite{Adams-Cockburn,Arnold-Awanou,Arnold-Awanou-Winther,
    Arnold-Winther-conforming, Awanou, Chen-Wang}; see
   \cite{Arnold-Winther-n, Gopalakrishnan-Guzman-n,HuShi,
    Man-Hu-Shi, Xie, Yi-3D, Yi} for nonconforming mixed finite elements,  and
   \cite{Arnold-Falk-Winther, Boffi-Brezzi-Fortin, Cockburn-Gopalakrishnan-Guzman,
    Gopalakrishnan-Guzman, Guzman} for new weakly symmetric finite elements.
However, most of these elements are difficult to be implemented;
    numerical examples can
    only be found in \cite{Carstensen-Eigel-Gedicke,
    Carstensen-Gunther-Reininghaus-Thiele2008,
   HuManZhang2012, HuManZhang2013, Yi} so far.

In a recent paper \cite{HuManZhang2013}, a new family of simple,
  any space-dimensional, symmetric, conforming mixed finite
  elements for the problem  is proposed.
In these elements,  quadratic polynomials $\{1, x_i, x_i^2\}$
   are used for the normal stresses $\sigma_{ii}$,
      bilinear polynomials $\{1, x_i, x_j, x_ix_j\}$
        for the shear stresses $\sigma_{ij}$,
        and  linear polynomials $\{1, x_i \}$ for
      the displacements $u_i$.
 The stress and displacement spaces of \cite{HuManZhang2013} are
  actually enrichment of those in \cite{HuManZhang2012},
   a family of symmetric nonconforming elements
These are possibly the simplest conforming mixed
  finite element methods.
A first order convergence was established
   for these elements in \cite{HuManZhang2013}.
However, superconvergence
  was observed from numerical examples presented  therein.

Superconvergence  is  one of the most active research fields
   for finite element methods.
A lot of  fundamental  results can be found  for  conforming,
   nonconforming and mixed finite elements
   of model problems in literature, see for instance,  \cite{Duran90,LinLin,LinYan}.
However,  no results can be found  for the mixed finite element methods
    under consideration in literature so far.
 A very recent paper \cite{ShiDongyang14}  analyzed  superconvergence
   of a family of conforming rectangular mixed  finite element methods
     for the two dimensional linear elasticity problem.
However, in the  conclusion,  it was pointed out that the technique therein
      can not be applied to mixed elements under consideration.

The aim of this paper is to prove superconvergence observed
      in \cite{HuManZhang2013}.
 One challenge is that
      the canonical interpolation operators for the stress spaces
     have no commuting properties,
 which are  indispensable  ingredients  for superconvergence analysis
     for mixed finite elements for the Poisson equations,
  see for instance, \cite{Brandts,DouglasWang1989,Duran90,EwingLazarovWang1991,Wang},
      also for the linear elasticity problem \cite{ShiDongyang14},
   and more details in Section 3.1.
Another challenge is that the normal stresses are coupled and
   consequently the superclose analysis  used in \cite{Duran90}
     for the mixed finite element of the Poisson equation can not be extended to
       the present case (see more details in Section 3.1.)
     To overcome these difficulties, we propose a new interpolation operator.
 Compared with the original interpolation
   operator from \cite{HuManZhang2013},  the new  one has a superclose
    property that is accomplished by adopting a new expansion
    which is motivated by a recent paper \cite{HuShi13}.
 Finally we propose an iterative  argument to establish an
     $O(h^{1+1/2})$ superconvergence.

This paper denotes by $H^k(T,X)$ the Sobolev space consisting of
   functions on domain $T\subset\mathbb{R}^n$, taking values in the
   finite-dimensional vector space $X$ and having all derivatives of
     order at most $k$ square-integrable.
For our purposes, the range
   space $X$ will be  $\mathbb{S}$, or $\mathbb{R}^n$, or
   $\mathbb{R}$.
In the latter case we may write simply $H^k(T)$.
$\|\cdot\|_{k,T}$ is the Sololev norm on $H^k(T)$.
Here $\mathbb{S}$ denotes
   the space of symmetric tensors, $H({\rm div},T,\mathbb{S})$,
   consisting of square-integrable symmetric matrix fields with
   square-integrable divergence.
The norm $\|\cdot\|_{H({\rm div},T)}$ reads
$$\|\tau\|_{H({\rm div},T)}^2:=\|\tau\|_{0,T}^2+\|{\rm div}\tau\|_{0,T}^2.$$
$L^2(T,\mathbb{R}^n)$ is the space of vector-valued functions
   which are square-integrable.

\section{The linear elasticity   problem and mixed finite elements}

\subsection{The linear elasticity   problem}
Based on the Hellinger-Reissner principle, the
  $n$-dimensional linear elasticity
  problem within a stress-displacement ($\sigma$-$u$) form reads:
Find $(\sigma,u)\in\Sigma\times V :=H({\rm div},\Omega,\mathbb
{S})
        \times L^2(\Omega,\mathbb{R}^n)$, such that
\an{\left\{ \ad{
  (A\sigma,\tau)+({\rm div}\tau,u)&= 0 && \hbox{for all \ } \tau\in\Sigma,\\
   ({\rm div}\sigma,v)&= (f,v) &\qquad& \hbox{for all \ } v\in V. }
   \right.\lab{eqn1}
}
Here the symmetric tensor space for stress $\Sigma$ and the
   space for vector displacement $V$ are, respectively,
  \a{
  H({\rm div},\Omega,\mathbb {S})
    &:= \Big\{ \p{\sigma_{ij} }_{n \times n} \in H(\d, \Omega)
    \ \Big| \ \sigma_{ij} = \sigma_{ji} \Big\}, \\
     L^2(\Omega,\mathbb{R}^n) &:=
     \Big\{ \p{u_1 & \cdots & u_n}^T
          \ \Big| \ u_i \in L^2(\Omega) \Big\}  .}

The matrix $A$ is defined as
   \a{
      A \sigma &= \frac 1{2\mu} \left(
       \sigma - \frac{\lambda}{2\mu + n \lambda} \operatorname{tr}(\sigma)
        \delta \right) }
  where $\delta$ is the identity matrix of $n\times n$, and $\mu$ and $\lambda$ are the
    Lam\'e constants.

This paper deals with a pure displacement problem with the
  homogeneous boundary condition that $u\equiv 0$ on
  $\partial\Omega$.
The domain is assumed to be a rectangular polyhedron
  in $\mathbb{R}^n$.

\subsection{The n-dimensional conforming mixed finite element space}

We recall  a conforming mixed finite element method proposed in \cite{HuManZhang2013}
 for the problem \meqref{eqn1}. We shall follow the notations used therein.

The rectangular domain $\Omega$ is subdivided by a family of
  rectangular grids  $\mathcal{T}_h$ (with the grid size $h$).
For convenience, the set of all $n-1$ dimensional faces in
  $\mathcal{T}_h$ is denoted by $\mathcal{F}_h$. For all element
  $K\in\mathcal{T}_h$, the set of all $n-1$ dimensional faces of $K$ perpendicular to $x_i$-axis
is denoted by $\mathcal{F}_{x_i,K}$, the set of all $n-2$
dimensional faces of $K$ perpendicular to $x_i$ and $x_j$ axes by
$\mathcal{F}_{x_i,x_j,K}$. Given any face
  $\mathcal{F}\in\mathcal{F}_h$, one fixed unit normal vector $\nu$
  with components $(\nu_1,\nu_2,\cdots,\nu_n)$ is assigned.


We first  introduce the finite element space locally on a single $n$-rectangle $
   K\in\mathcal{T}_h$:
  \a{
   V(K):=\Big\{ v =(v_1,\dots, v_n) \ \Big| \
       v_i \in P_1(x_i) \Big\},}
   \a{
        \Sigma(K):  =\Big\{ \sigma\in\p{ \sigma_{ij} }_\mathbb{S}
   \  \Big| \ \sigma_{ii}\in P_2(x_i),
    \sigma_{ij}\in Q_1(x_i,x_j),i\neq j \Big\}, }
where $$\begin{array}{rcl} P_1(x_i)&:=&{\rm span}\left\{1,x_i\right\},\\
P_2(x_i)&:=&{\rm span}\left\{1,x_i,x_i^2\right\},\\
Q_1\left(x_i,x_j\right)&:=&{\rm
span}\left\{1,x_i,x_j,x_ix_j\right\}.
\end{array}$$ For example, in 2D ($n=2$), the spaces may be
displayed as
 \a{ \Sigma(K)&=\operatorname{span}
    \p{ \left\{1, x_1, x_1^2 \right\}
            & \left\{1, x_1, x_2, x_1x_2 \right\}
             \\
        \left\{1, x_1, x_2, x_1x_2 \right\}
            & \left\{1, x_2, x_2^2 \right\}
             }, \\
      V(K)& =\operatorname{span}
    \p{ \left\{1, x_1 \right\}  \\
        \left\{1, x_2 \right\}  }.
  }

Due to the $H(\d)$ requirement, $\sigma_{ii}$ has to be continuous
in $x_i$ direction,
   while $\sigma_{ij}$ has to be continuous in both $x_i$ and $x_j$ directions.
Thus, we can specify the local degrees of freedom for the two
finite element spaces on element $K$ as follows,
\a{
 \bullet & \quad \disp\frac{1}{|K|}\int_K u_i~v{\rm d}V,&&\hbox{for all \ }
    v\in P_1(x_i)  & \hbox{ and \ }& i=1,2,\dots, n; && \\
  \bullet & \quad \disp\frac{1}{|F_{x_i, K}|}\int_{F_{x_i, K}}\sigma_{ii}{\rm
          d}F,&&\hbox{for all \ }
          F_{x_i,K}\in\mathcal{F}_{x_i,K}  & \hbox{ and \ }&
            i=1,2,\cdots, n;&& \\
  \bullet & \quad\frac{1}{|K|} \int_K\sigma_{ii}{\rm d}V,&&  &&   i=1,2,\cdots, n;&& \\
  \bullet & \quad\frac{1}{|F_{x_i,x_j,K}|} \int_{F_{x_i,x_j,K}}\sigma _{ij} {\rm
        d}F,&&\hbox{for all \ }
         F_{x_i,x_j,K}\in\mathcal{F}_{x_i,x_j,K}  & \hbox{ and \ }&
           1\leq i<j\leq n.
    }

The global spaces $\Sigma_h$ and $V_h$ can be defined by their
   property
   \an{\lab{Sh}
 \Sigma_h&:=\{~\sigma
    \in H(\d,\Omega,\mathbb{S}) \ | \ \sigma|_K\in\Sigma(K)
            \ \hbox{ for all } K\in\mathcal{T}_h  \}, \\
  \lab{Vh}
   V_h &:= \{v\in L^2(\Omega,\mathbb{R}^n) \ | \
         v|_K\in V(K) \ \hbox{ for all } K\in\mathcal{T}_h \}.
    }

The mixed finite element approximation of Problem (\ref{eqn1})
reads: Find
   $(\sigma_h,~u_h)\in\Sigma_h\times V_h$ such that
 \e{ \left\{ \ad{
    (A\sigma_h, \tau)+({\rm div}\tau, u_h)&= 0 &&
              \hbox{for all \ } \tau \in\Sigma_h,\\
     (\d\sigma_h, v)& = (f, v) &&  \hbox{for all \ } v\in V_h.
      } \right. \lab{DP}
    }

It follows from the definition of $V_h$ and $\Sigma_h$ that
   \a{ \d  \Sigma_h \subset V_h.}
This, in turn, leads to a strong divergence-free space:
 \an{ \lab {kernel}
    Z_h&:= \{\tau_h\in\Sigma_h \ | \ (\d\tau_h, v)=0 \quad
    \hbox{for all } v\in V_h\}\\
    \nonumber
          &= \{\tau_h \in\Sigma_h \ | \  \d \tau_h=0
    \hbox{\ pointwise } \}.
    }
\subsection{Well-posedness of the discrete problem}
The well-posedness of the discrete problem
    \meqref{DP} is proved in \cite{HuManZhang2013}. More precisely, it was shown therein that:

\begin{enumerate}
\item K-ellipticity. There exists a constant $C>0$, independent of the
   meshsize $h$, such that
    \an{ \lab{below} (A\tau, \tau)\geq C\|\tau\|_{H(\d)}^2\quad
       \hbox{for all } \tau \in Z_h, }
    where $Z_h$ is the divergence-free space defined in \meqref{kernel}.

\item  Discrete B-B condition.
    There exists a positive constant $C>0$
            independent of the meshsize $h$, such that
    \an{\lab{inf-sup}
   \inf_{0\neq v\in V_h}   \sup_{0\neq\tau\in\Sigma_h}\frac{({\rm
        div}\tau, v)}{\|\tau\|_{H(\d)}  \|v\|_{0} }\geq
    C .}
\end{enumerate}

In addition, there is a refined  discrete inf--sup condition from \cite{HuManZhang2013} as follows:
\begin{lemma}\lab{lemma1}
For any $v\in V_h$, there exists  an $H(\d)$ field
\a{ \tau =\p{ \tau_{11}  &  0&\cdots&\cdots&0\\
                 0   & \ddots& & &\vdots\\
                \vdots&    &\tau_{ii}& &\vdots \\
                \vdots& & &\ddots& 0 \\
                0&\cdots&\cdots&0&\tau_{nn}}\in \Sigma_h,
                  }
such that
  \an{\lab{equal} \d \tau=v \text{ and }\frac{({\rm
        div}\tau, v)}{\|\tau\|_{H(\d)}}\geq \sqrt{\frac 23}\|v\|_0.
   }
\end{lemma}

\subsection{Error estimate} Since it is very  difficult to show the superclose property of the interpolation operator
 defined in \cite{HuManZhang2013},  for any $\sigma\in H^2(\Omega,\mathbb{S})$, we define a new  interpolation by
\a{
 \Pi_h\sigma
   =\p{ \Pi_{11}\sigma_{11}  & \Pi_{12}\sigma_{12}
   &\cdots&\cdots&\Pi_{1n}\sigma_{1n}\\
                \vdots & \ddots&\vdots&\vdots&\vdots\\
                \Pi_{i1}\sigma_{i1}&\cdots&\Pi_{ii}\sigma_{ii}
    &\cdots&\Pi_{in}\sigma_{in} \\
                \vdots & \vdots&\vdots&\ddots&\vdots\\
                \Pi_{n1}\sigma_{n1}&\Pi_{n2}\sigma_{n2}
    &\cdots&\cdots&\Pi_{nn}\sigma_{nn}}\in\Sigma_h, }
   where $\Pi_{ij}=\Pi_{ji}$, are defined next.
The interpolation operator $\Pi_{ii}$ is defined by, for any
     $K\in\mathcal{T}_h$,
 \a{
    \int_{F_{x_i,K}}\Pi_{ii}\sigma_{ii}dF & =\int_{F_{x_i,K}}\sigma_{ii}dF
        \quad  &&\hbox{for all \ }
   F_{x_i,K}\in\mathcal{F}_{x_i,K},\\
   \int_K\Pi_{ii}\sigma_{ii}dV & =\int_K\sigma_{ii}dV && \hbox{for all \ }
    ~K\in\mathcal{T}_h. }
    Here and  throughout this paper, $\mathcal{F}_{x_i,K}$ denotes the set of $n-1$ dimensional faces of $K$ which are perpendicular
     to the $x_i$ axis and $\mathcal{F}_{x_i}=\cup_{K\in\mathcal{T}_h}\mathcal{F}_{x_i,K}$.

To define $\Pi_{ij}$, we introduce the nodal basis functions on unit element.
  \a{
   \phi_{k}(x,y) :=
    \begin{cases}\ \ \, (x-1)(y-1), & k=0, \\
      -(x-0)(y-1), & k=1, \\
       \ \ \,    (x-0)(y-0), & k=2, \\
      -(x-1)(y-0), &k=3. \end{cases} \\
   }
We also use multi-index notations as follows
 \an{\lab{out2}   \sum_{L_{\text{\sout{$i$}} ,\text{\sout{$j$}} }=1}^N :=
    \underset {\{1\le l_k \le N \mid k\ne i, k\ne j\}  }
     {\sum\cdots\sum }. }

Since $\sigma_{ij}\in H^2(\Omega)$, these  basis functions and the short
   notation allow for  defining the interpolation as follows
    \begin{equation}
   \Pi_{ij} \sigma_{ij}  =
\sum_{L_{\text{\sout{$i$}} , \text{\sout{$j$}} }=1}^N
         \sum_{ l_i,l_j=1}^N  \sum_{k=0}^3
       c^{(ij),k}_{l_1,\cdots,l_n} \phi_k(\frac{x_i}h-(l_i-1) ,
        \frac{x_j}h -(l_j-1) )
    \end{equation}
 where the interpolation parameters satisfy
 \a{ &\quad \ c^{(ij),2}_{l_1,\cdots,(l_i-1),\cdots,(l_j-1),\cdots,l_n}\\
     &= c^{(ij),3}_{l_1,\cdots,l_i,\cdots,(l_j-1),\cdots, l_n}
     = c^{(ij),0}_{l_1,\cdots,l_i,\cdots,l_j,\cdots, l_n}
     = c^{(ij),1}_{l_1,\cdots,(l_i-1),\cdots,l_j,\cdots, l_n}\\
    &=\frac{1}{|F_{x_i, x_j, l_1,\cdots, l_i, \cdots, l_j, \cdots, l_n}|}
     \int_{F_{x_i, x_j, l_1,\cdots, l_i, \cdots, l_j, \cdots, l_n}}\sigma_{ij}dF,
       \qquad 0<l_i,l_j<N, }
where $F_{x_i, x_j, l_1,\cdots, l_i, \cdots, l_j, \cdots, l_n}$ is the unique
      $n-2$ dimensional face at vertex
$$((l_1-1)h,\cdots, (l_i-1)h, \cdots, (l_j-1), \cdots, (l_n-1)h)$$
 which is shared by elements:
 \begin{equation*}
 \begin{split}
 K_1&=[l_1,\cdots, l_i, \cdots, l_j, \cdots, l_n],\\
  K_2&=[l_1,\cdots,(l_i-1),\cdots,(l_j-1),\cdots,l_n],\\
  K_3&=[l_1,\cdots,l_i,\cdots,(l_j-1),\cdots, l_n],\\
  K_4&=[l_1,\cdots,(l_i-1),\cdots,l_j,\cdots, l_n],\\
  \end{split}
  \end{equation*}
   $|F_{x_i, x_j, l_1,\cdots, l_i, \cdots, l_j, \cdots, l_n}|$ is the measure
 of face $F_{x_i, x_j, l_1,\cdots, l_i, \cdots, l_j, \cdots, l_n}$.  If
 $$|F_{x_i, x_j, l_1,\cdots, l_i, \cdots, l_j, \cdots, l_n}|=0,$$
  $$\frac{1}{|F_{x_i, x_j, l_1,\cdots, l_i, \cdots, l_j, \cdots, l_n}|}
    \int_{F_{x_i, x_j, l_1,\cdots, l_i, \cdots, l_j, \cdots, l_n}}\sigma_{ij}dF$$
   is understood as the value of $\sigma_{ij}$ at  vertex
   $$
   ((l_1-1)h,\cdots, (l_i-1)h, \cdots, (l_j-1), \cdots, (l_n-1)h).
   $$ The operator $\Pi_{ij}$
    is different from that defined in \cite{HuManZhang2013}.
Note that the corresponding operator of \cite{HuManZhang2013} can not be used
     for the superconvergence analysis.
For this operator, we have the following error estimates:
    \an{\lab{inter-error1}
  \|\sigma_{ij}-\Pi_{ij}\sigma_{ij}\|_{0,K} & \leq Ch\|\sigma_{ij}\|_{2,K} ,\\
   \lab{inter-error2}
\|\frac{\partial}{\partial x_i}(\sigma_{ij}-\Pi_{ij}\sigma_{ij})\|_{0,K}
   & \leq Ch\|\sigma_{ij}\|_{2,K} , \\
  \lab{inter-error3}
  \|\frac{\partial}{\partial x_j}(\sigma_{ij}-\Pi_{ij}\sigma_{ij})\|_{0,K}
   &\leq Ch\|\sigma_{ij}\|_{2, K}.}

Since the space for the operator $\Pi_{ii}$ contains
 the 1D quadratic polynomials ${\rm span}\{1,x_i,x_i^2\}$, the
 scaling argument and standard approximation,  state
  \an{\lab{IE1}
   |\sigma_{ii}-\Pi_{ii}\sigma_{ii}|_{0,K}  & \leq
    Ch|\sigma_{ii}|_{1,K} , \\
    \lab{IE2} |\frac{\partial}{\partial
    x_i}(\sigma_{ii}-\Pi_{ii}\sigma_{ii})|_{0,K} & \leq
    Ch|\frac{\partial{\sigma_{ii}}}{\partial
     x_i}|_{1,K}, }
 for all $ K\in\mathcal{T}_h$.

A summary of these aforementioned estimates
(\ref{inter-error1})--(\ref{IE2}) leads to
\begin{theorem}\lab{theorem-inter}
   For any $\sigma\in H^2(\Omega, \mathbb{S})$, we have that
 \an{ \lab{inter-err4}
   \|\sigma-\Pi_h\sigma\|_0 & \leq Ch\|\sigma\|_2, \\
  \lab{inter-err5}
   \|{\rm
   div}(\sigma-\Pi_h\sigma)\|_0& \leq Ch\|\sigma\|_2. }
   \end{theorem}

The stability of the elements and the standard theory of mixed
  finite element methods, see for instance \cite{Brezzi, Brezzi-Fortin}, give the
  following abstract error estimate: \an{
  \lab{theorem-err1} \|\sigma-\sigma_h\|_{H({\rm
  div})}+\|u-u_h\|_0\leq C \inf\limits_{\tau_h\in\Sigma_h,v_h\in
  V_h}\left(\|\sigma-\tau_h\|_{H({\rm div})}+\|u-v_h\|_0\right).}
Let $P_h$ denotes the projection operator from $V$ to $V_h$,
  which has the error estimate
\an{\lab{proj-error}
   \|v-P_hv\|_0\leq Ch\|v\|_1. }
Choosing $\tau_h=\Pi_h\sigma$ and
  $v_h=P_hu$ in (\ref{theorem-err1}), the estimates
  (\ref{inter-err4}), (\ref{inter-err5}), (\ref{proj-error}) prove
\begin{theorem}\lab{MainError} Let
  $(\sigma, u)\in\Sigma\times V$ be the exact solution of
   problem \meqref{eqn1} and $(\tau_h, u_h)\in\Sigma_h\times
   V_h$ the finite element solution of \meqref{DP}.  Then,
\a{
\|\sigma-\sigma_h\|_{H({\rm div})}
    &\le Ch(\|\sigma\|_2+\|u\|_1),\\
 \|u-u_h\|_0&\le     Ch(\|\sigma\|_2+\|u\|_1).
      }
\end{theorem}

\section{The superclose property of the canonical interpolations}
\subsection{Main difficulties}
By the K--ellipticity  in \eqref{below} and the discrete
   inf--sup condition \eqref{inf-sup}, it is routine to prove that
\begin{equation}\lab{eq3.1}
\begin{split}
&\quad \ \|\Pi_h\sigma-\sigma_h\|_{H(\d)}+\|P_hu-u_h\|_0\\
&\leq C\sup\limits_{0\ne(\tau, v)\in\Sigma_h\times V_h}
   \frac{ (A(\sigma-\Pi_h\sigma), \tau)+(u-P_h u, \d \tau)
     -(\d (\sigma-\Pi_h\sigma), v) }{\|\tau\|_{H(\d)}+\|v\|_0}.
\end{split}
\end{equation}
Note that the inequality \eqref{eq3.1} is the starting point for
    superconvergence analysis of mixed finite element methods, see
 for instance, \cite{Brandts,Duran90} and \cite{ShiDongyang14}.
However, this formulation can not be directly used for mixed finite
    elements under consideration,  the reasons lie in that
  \begin{itemize}
  \item The interpolation operator lacks the usual commuting property, namely,
   $$
   \d \Pi_h\sigma \not =P_h\d \sigma.
   $$
   \item The components of the stress normal are  coupled  through  $A(\sigma-\Pi_h\sigma)$, it is impossible  to prove directly
    the following super-close property
   $$
   (A(\sigma-\Pi_h\sigma), \tau)\leq Ch^2\|\sigma\|_2\|\tau\|_{H(\d)}
   $$
   for a general $\tau\in\Sigma_h$.
  \end{itemize}
  In the sequel, we will need some results on Sobolev spaces. They are formulated
in the following lemma. First of all, define $\partial\Omega_h$ as the subset of points
 having (Euclidian) distance less than $h$ from the boundary:
 $$
 \partial\Omega_h:=\{x\in\Omega|  \exists y\in\partial\Omega: {\rm dist}(x,y)\leq h \}.
 $$
  \begin{lemma}\lab{trace}
  For $v\in H^s(\Omega)$ with $0\leq s\leq 1/2$, it holds
  \begin{equation}
  \|v\|_{0,\partial\Omega_h}\leq Ch^s\|v\|_s.
  \end{equation}
  \end{lemma}

  \subsection{The superclose property of  $(\d (\sigma-\Pi_h\sigma), v)$}
To overcome the first difficulty, we  follow the idea of \cite{HuShi13}
   to adopt a new expansion of the operator $\Pi_h$.
In fact,
 let $\Pi_K=\Pi_h|_K$, we have the following crucial result.
 \begin{lemma} For any $\sigma\in P_2(K, \mathbb{S})$ and $v\in V_K$, it holds that
 \begin{equation}\lab{eq3.2}
 (\d (\sigma-\Pi_K\sigma), v)_K=0.
 \end{equation}
 \end{lemma}
 \begin{proof} We only need to prove the result on the reference element $K=[-1,1]^n$.
For any $\sigma\in P_2(K, \mathbb{S})$, its components $\sigma_{i j}$, $i, j=1, \cdots, n$
    with $i\not= j$, can be expressed as
 \begin{equation*}
 \sigma_{ij}=p_0(x_i, x_j)+x_ip_1^{(ij)}+x_jp_2^{(ij)}+p_3^{(ij)},
 \end{equation*}
 where $p_0(x_i, x_j)$ is a polynomial of degree $2$ in $x_i$ and $x_j$,
   both $p_1^{(ij)}$ and $p_2^{(ij)}$ are homogeneous
    polynomials of degree $1$ of $(n-2)$ of variables $x_k$, $k=1, \cdots, n$ with $k\not=i, j$,
   and $p_3^{(ij)}\in P_2(K)$
     is a polynomial of degree $2$ with respect to $(n-2)$
        variables $x_k$, $k=1, \cdots, n$ with $k\not=i, j$.
The definition of $\Pi_{ij, K}$ leads to
    \begin{equation}\lab{eq3.3}
    \begin{split}
    \sigma_{ij}-\Pi_{ij, K}\sigma_{ij}&=c_{ii}(x_i^2-1)+c_{ij}
              (x_j^2-1)+x_ip_1^{(ij)}+x_jp_2^{(ij)}\\
    &\quad +p_3^{(ij)}-\frac{1}{|F_{x_i, x_j}|}\int_{F_{x_i, x_j}}p_3^{(ij)}dF,
   \end{split}
    \end{equation}
     for two interpolation parameters $c_{ii}$ and $c_{ij}$, where $F_{x_i, x_j}$ is any
    $n-2$ dimensional face  of $K$ which is perpendicular  to the plane $\sspan \{ x_i, x_j\}$.
Here we use the facts that $\int_{F_{x_i, x_j}}p_1^{(ij)}dF=\int_{F_{x_i, x_j}}p_2^{(ij)}dF=0$,
  and that $p_3^{(ij)}$ is a constant function with respect to variables
     $x_i$ and $x_j$, and that $\sum\limits_{k=0}^3\phi_{k}=1$.

    If $i=j$, we have
 \begin{equation*}
 \sigma_{ii}=p_0(x_i)+x_ip_1^{(ii)}+p_2^{(ii)},
 \end{equation*}
 where $p_0(x_i)$ is a polynomial of degree $2$ or less in one variable $x_i$,
   and $p_1^{(ii)}$ is
  a homogeneous polynomial of degree $1$ of $(n-1)$  variables $x_k$,
    $k=1, \cdots, n$ with $k\not=i$, and $p_2^{(ii)}\in P_2(K)$
   is a polynomial of degree $2$ of $(n-1)$ variables $x_k$,
     $k=1, \cdots, n$ with $k\not=i$.
The definition
    of $\Pi_{ii, K}$ leads to
 \begin{equation}\lab{eq3.4}
 \sigma_{ii}-\Pi_{ii, K}\sigma_{ii}=x_i p_1^{(ii)}+p_2^{(ii)}-\frac{1}{|K|}\int_Kp_2^{(ii)}dV.
 \end{equation}

  By \eqref{eq3.3} and \eqref{eq3.4}, we have
   \begin{equation*}
        \frac{\partial( \sigma_{ij}-\Pi_{ij, K}\sigma_{ij})}{\partial x_j}=p_2^{(ij)}+2c_{ij}x_j
      \hbox{ \     and \ }
        \frac{\partial(\sigma_{ii}-\Pi_{ii, K}\sigma_{ii})}{\partial x_i}=p_1^{(ii)}.
 \end{equation*}
 Hence,  the $i$-th component of $\d (\sigma-\Pi_K\sigma)$ can be expressed as
  \begin{equation*}
 (\d (\sigma-\Pi_K\sigma))_i=p_1^{(ii)}+\sum\limits_{i\not=j=1}^n
     ( p_2^{(ij)}+2 c_{ij}x_j).
  \end{equation*}
We can compute
\a{ \int_K (\d (\sigma-\Pi_K\sigma))_i dV &= \int_K p_1^{(ii)} dV+\sum\limits_{i\not=j=1}^n
     \int_K( p_2^{(ij)}   +2 c_{ij}x_j ) dV\\
               &= 0 + \sum\limits_{i\not=j=1}^n
     ( \int_Kp_2^{(ij)} dV +2 c_{ij}\int_Kx_j dV) = 0, \\
    \int_K (\d (\sigma-\Pi_K\sigma))_i x_i dV &= \int_{-1}^1 x_i d x_i
           \int_{K_{n-1}} \bigg( p_1^{(ii)}+\sum\limits_{i\not=j=1}^n
     ( p_2^{(ij)}  +2 c_{ij}x_j  ) \bigg) d V_{n-1} \\
               &= 0  \cdot 0  = 0, }
   Here $K_{n-1}=[0,1]^{n-1}$ is the $n-1$ dimensional cube without variable $x_i$.
  Note that the $i$-th component of $v$ can be written as
  \begin{equation*}
  v_i=a_0+a_1x_i
  \end{equation*}
  for two interpolation parameters $a_0$ and $a_1$.
Thus
$$
  ((\d (\sigma-\Pi_K\sigma))_i, v_i)_K=0,
$$
which completes the proof.
 \end{proof}

 As a consequence of \eqref{eq3.2},  we have the following
     superclose property for the term $(\d (\sigma-\Pi_h\sigma), v)$.
\begin{lemma} Suppose that $\sigma\in H^3(\Omega, \mathbb{S})$. Then it holds that
\begin{equation}\lab{eq3.5}
|(\d (\sigma-\Pi_h\sigma), v)|\leq Ch^2|\sigma|_3\|v\|_0 \text{ for any }v\in V_h.
\end{equation}
\end{lemma}
\begin{proof} Given element $K$, let $I_{2,K}: H(\d, K, \mathbb{S})\rightarrow
    P_2(K, \mathbb{S})$ be the $L^2$ projection operator defined
    as: Given $\tau\in H(\d, K, \mathbb{S})$,
  find $I_{2,K}\tau\in P_2(K, \mathbb{S})$ such that
 \begin{equation*}
 \int_K I_{2,K}\tau q dV=\int_K\tau qdV\text{ for any }q\in P_2(K, \mathbb{S}).
 \end{equation*}
This allows for the following decomposition:
\a{ (\d (\sigma-\Pi_h\sigma), v) & =\sum\limits_{K\in\mathcal{T}_h}
      (\d ((I-\Pi_K)I_{2,K}\sigma+(I-\Pi_K)(I-I_{2,K})\sigma), v)_K\\
      & =\sum\limits_{K\in\mathcal{T}_h}
      (\d ((I-\Pi_K)(I-I_{2,K})\sigma), v)_K, }
   where we applied \eqref{eq3.2}.
    The desired result follows from the stability of $\Pi_h$
        and the approximation property of $I_{2,K}$.
\end{proof}

\subsection{The superclose property of $(A(\sigma-\Pi_h\sigma), \sigma_h-\Pi_h\sigma)$}

  To deal with the second difficulty, we propose to
explore  the strong discrete inf--sup condition presented in Lemma \ref{lemma1}.

\begin{lemma}\lab{lemma3.3} For any $\sigma\in P_1(K, \mathbb{S})$
    and $\tau\in \Sigma_{n,h}$, it holds that
 \begin{equation}\lab{eq3.6}
     |(A (\sigma-\Pi_K\sigma), \tau)_K |
           \leq Ch^2 \|\d\tau\|_{0, K}|\sigma|_{1, K},
 \end{equation}
 where
 $$
     \Sigma_{n,h}=\{\tau=\diag(\tau_{11}, \cdots, \tau_{nn}), \
             \tau\in \Sigma_{h}\}.
 $$
\end{lemma}
\begin{proof} We only need to prove the result on the reference element
    $K=[-1,1]^n$. For any $\sigma\in P_1(K, \mathbb{S})$, its normal
      components can be written as
 \begin{equation*}
    \sigma_{ii}=c_0^{(ii)}+\sum\limits_{j=1}^nc_j^{(ii)}x_j, \quad i=1, \cdots, n,
 \end{equation*}
 where $c_j^{(ii)}$, $j=0, \cdots, n$, are interpolation parameters.
By the definition of the operator $A$,  the $ii$-th component
  of $A(\sigma-\Pi_K\sigma)$ is
  $$
  A(\sigma-\Pi_K\sigma)_{ii}=\frac{1}{2\mu(2\mu+n\lambda)}\bigg((2\mu+(n-1)\lambda)
  \sum\limits_{i\not=j=1}^nc_j^{(ii)}x_j
 -\lambda\sum\limits_{i\not=k=1}^n\sum\limits_{k\not=j=1}^nc_j^{(kk)}x_j\bigg).
  $$
Note that the $i$-th component of $\tau$ can be written as
  \begin{equation*}
  \tau_{ii}=a_0^{(ii)}+a_1^{(ii)}x_i+a_2^{(ii)}x_i^2
  \end{equation*}
  for  parameters $a_0^{(ii)}$,  $a_1^{(ii)}$ and $a_2^{(ii)}$. Therefore,
  \an{\lab{eq3.7-2}
    \ad{(A(\sigma-\Pi_K\sigma)_{ii}, \tau_{ii})_K
   &=-a_1^{(ii)}\frac{\lambda}{2\mu(2\mu+n\lambda)}\sum\limits_{i\not=k=1}^nc_i^{(kk)}(x_i,x_i)\\
  &=-\frac{2^n}{3}\frac{\partial\tau_{ii}}{\partial x_i}(0)
     \frac{\lambda}{2\mu(2\mu+n\lambda)}\sum\limits_{i\not=k=1}^n
     \frac{\partial \sigma_{kk}}{\partial x_i},} }
  \a{ |(A(\sigma-\Pi_K\sigma)_{ii}, \tau_{ii})_K|
    &\leq C \frac{\lambda}{6\mu(2\mu+n\lambda)}\|\frac{\partial\tau_{ii}}{\partial x_i}\|_{0, K}
     \|\sum\limits_{i\not=k=1}^n\frac{\partial \sigma_{kk}}{\partial x_i}\|_{0, K}. }
A summation over all $n$ components leads to
  \begin{equation*}
  \begin{split}
   |(A(\sigma-\Pi_K\sigma), \tau)_K|\leq C\|\d\tau\|_{0,K}|\sigma|_{1,K}.
  \end{split}
  \end{equation*}
  The final result follows from a scaling argument.
\end{proof}

A combination of  the above lemma and \eqref{eq3.5} yields  the following important result.

\begin{lemma} It holds that
\begin{equation}\lab{eq3.7-1}
\|\Pi_h\sigma-\sigma_h\|_0^2+\|P_hu-u_h\|_0^2
    \leq C(A(\Pi_h\sigma-\sigma), \Pi_h\sigma-\sigma_h) +Ch^4\|\sigma\|_3^2.
\end{equation}
\end{lemma}
\begin{proof}
It follows from Lemma \ref{lemma1} that there exists
     $\tau=\diag(\tau_{11}, \cdots, \tau_{nn})\in\Sigma_h$ such that
\begin{equation}\lab{eq3.7}
\d \tau=u_h-P_hu \text{\ and \ }\|\tau\|_{H(\d)}\leq C\|u_h-P_hu\|_0.
\end{equation}
This allows for the following decomposition:
\begin{equation*}
\begin{split}
(u_h-P_hu, u_h-P_hu)&=(u_h-P_hu, \d \tau)=(u_h-u, \d\tau)=(A(\sigma-\sigma_h), \tau)\\
&=(A(\sigma-\Pi_h\sigma), \tau)+(A(\Pi_h\sigma-\sigma_h), \tau).
\end{split}
\end{equation*}
Since $\tau$ is a diagonal matrix, it follows from \eqref{eq3.6}  and \eqref{eq3.7} that
\begin{equation*}
(A(\sigma-\Pi_h\sigma), \tau)\leq Ch^2|\sigma|_1\|\d\tau\|_0\leq Ch^2|\sigma|_1\|u_h-P_hu\|_0.
\end{equation*}
A substitution of this inequality into the previous equation,
   by the Cauchy--Schwarz inequality and \eqref{eq3.7}, leads to
\begin{equation}\lab{eq3.8}
\|u_h-P_hu\|_0\leq C(h^2|\sigma|_1+\|\Pi_h\sigma-\sigma_h\|_0).
\end{equation}
On the other hand, let $\tau=\Pi_h\sigma-\sigma_h$, we have
\begin{equation}
\begin{split}
(A(\Pi_h\sigma-\sigma_h), \tau)&=(A(\Pi_h\sigma-\sigma), \tau)+(A(\sigma-\sigma_h), \tau)\\
&=(A(\Pi_h\sigma-\sigma), \tau)-(u-u_h, \d  \tau)\\
&=(A(\Pi_h\sigma-\sigma), \tau)-(P_hu-u_h, \d  (\Pi_h\sigma-\sigma_h))\\
&=(A(\Pi_h\sigma-\sigma), \tau)-(P_hu-u_h, \d  (\Pi_h\sigma-\sigma)).
\end{split}
\end{equation}
From \eqref{eq3.5} and \eqref{eq3.8} it follows
\begin{equation}
\begin{split}
(A(\Pi_h\sigma-\sigma_h), \tau)&\leq (A(\Pi_h\sigma-\sigma), \tau)
    +Ch^2|\sigma|_3\|P_hu-u_h\|_0\\
&\leq (A(\Pi_h\sigma-\sigma), \tau)+Ch^2|\sigma|_3(h^2|\sigma|_1+\|\Pi_h\sigma-\sigma_h\|_0).
\end{split}
\end{equation}
Since there exists a positive constant $\beta$ such that
$$
\beta\|\Pi_h\sigma-\sigma_h\|_0^2\leq (A(\Pi_h\sigma-\sigma_h), \tau),
$$
an application of the Young inequality leads to
$$
\|\Pi_h\sigma-\sigma_h\|_0^2+\|P_hu-u_h\|_0^2
   \leq C(A(\Pi_h\sigma-\sigma), \tau) +Ch^4\|\sigma\|_3^2,
$$
which completes the proof.
\end{proof}

\begin{lemma} For any $\sigma_{ij}\in P_1(K)$ and $\tau_{ij}\in Q_1(x_i, x_j)$, it holds that
\begin{equation}
(\sigma_{ij}-\Pi_{ij, K}\sigma_{ij}, \tau_{ij})_K=0.
\end{equation}
\end{lemma}
\begin{proof}
We only need to prove the result on the reference element $K=[-1,1]^n$.
Since $\sigma_{ij}\in P_1(K)$ , we have
 \begin{equation*}
    \begin{split}
    \sigma_{ij}-\Pi_{ij, K}\sigma_{ij}&=p^{(ij)},
   \end{split}
    \end{equation*}
 where $p^{(ij)}$ be a homogeneous polynomial of degree $1$ with respect to
      variables $x_k$, $k=1, \cdots, n$ but $k\not=i, j$.
   Any $\tau_{ij}\in Q_1(x_i, x_j)$ can be expressed as
 \begin{equation*}
 \tau_{ij}=a_0+a_1x_i+a_2x_j+a_3x_ix_j,
 \end{equation*}
 for four interpolation parameters $a_k$, $k=0, \cdots, 3$.  On the reference element $K$, it is straightforward to see that
\begin{equation*}
 (\sigma_{ij}-\Pi_{ij, K}\sigma_{ij}, \tau_{ij})=0.
 \end{equation*}
 This completes the proof.
\end{proof}

This lemma and a similar argument of \eqref{eq3.5} can prove the following supercloseness.
\begin{lemma} For any $\tau_{ij}\in \Sigma_{ij,h}
    :=\{e_i^\prime\tau e_j, \tau\in\Sigma_h \}$ it holds that
\begin{equation}\lab{eq3.15}
(\sigma_{ij}-\Pi_{ij}\sigma_{ij}, \tau_{ij})\leq Ch^2|\sigma_{ij}|_2\|\tau_{ij}\|_0,
\end{equation}
provided that $\sigma_{ij}\in H^2(\Omega)$. Here $e_i$ and $e_j$ are the $i$-th and $j$-th  canonical basis of the space $\mathbb{R}^n$, respectively.
\end{lemma}

\begin{lemma}\lab{lemma3.8} Let $(\sigma, u)$ and $(\sigma_h, u_h)$
    be solutions of problems \eqref{eqn1} and \eqref{DP},
  respectively.
Suppose that $\sigma \in H^2(\Omega, \mathbb{S})$ and $u\in H^1(\Omega, \mathbb{R}^n)$.
  Then  there holds that
\begin{equation}\lab{eq3.16}
   (A(\sigma-\Pi_h\sigma), \sigma_h-\Pi_h\sigma)
    \leq Ch^{5/2}(\|\sigma\|_2+\|u\|_1)\|\sigma\|_2.
\end{equation}
\end{lemma}
\begin{proof} Let $\tau=\sigma_h-\Pi_h\sigma$. Given element $K$, let $I_{1,K}: L^2(K)\rightarrow
    P_1(K)$ be the $L^2$ projection operator defined
    as: Given $v \in L^2(K)$,
  find $I_{1,K}v\in P_1(K)$ such that
 \begin{equation*}
 \int_K I_{1,K}v q dV=\int_Kv qdV\text{ for any }q\in P_1(K, \mathbb{S}).
 \end{equation*}
This leads to the following decomposition:
\begin{equation*}
  \begin{split}
  &(A(\sigma-\Pi_h\sigma)_{ii}, \tau_{ii})
  =\sum\limits_{K\in\mathcal{T}_h}(A(\sigma-\Pi_K\sigma)_{ii}, \tau_{ii})_K\\
  &=\sum\limits_{K\in\mathcal{T}_h}(A((I-\Pi_K)I_{1,K}\sigma)_{ii}, \tau_{ii})_K+\sum\limits_{K\in\mathcal{T}_h}(A((I-\Pi_K)(I-I_{1,K})\sigma)_{ii}, \tau_{ii})_K.
  \end{split}
  \end{equation*}
Then it follows from \eqref{eq3.7-2} that
\begin{equation*}
  \begin{split}
  &(A(\sigma-\Pi_h\sigma)_{ii}, \tau_{ii})
  =-\frac{\lambda \, h^2}{24 \mu(2\mu+n\lambda)}\sum\limits_{i\not=k=1}^n
  \sum\limits_{K\in\mathcal{T}_h} \big(\frac{\partial (I_{1,K}\sigma)_{kk}}{\partial x_i}, \frac{\partial\tau_{ii}}{\partial x_i}\big)_K
    +Ch^2\|\sigma\|_2\|\tau_{ii}\|_0.
  \end{split}
  \end{equation*}
 Since $\frac{\partial\tau_{ii}}{\partial x_i}$ is of the form  $a_1^{(ii)}+a_2^{(ii)}x_i$  for  parameters  $a_1^{(ii)}$ and $a_2^{(ii)}$, $\big(\frac{\partial (I_{1,K}x_l^2-x_l^2)}{\partial x_i}, \frac{\partial\tau_{ii}}{\partial x_i}\big)_K=0$ for $l\not=i$.  Therefore,
     $$
  \big(\frac{\partial (I_{1,K}\sigma-\sigma)_{kk}}{\partial x_i}, \frac{\partial\tau_{ii}}{\partial x_i}\big)_K=-\frac{h}{3}\big(\frac{\partial^2 \sigma_{kk}}{\partial x_i^2}, \frac{\partial\tau_{ii}}{\partial x_i}\big)_K+Ch^2|\sigma|_{3,K}|\tau_{ii}|_{1,K}.
 $$
After an elementwise inverse estimate, a combination of these two equations yield
\begin{equation}\lab{eq3.18}
  \begin{split}
  (A(\sigma-\Pi_h\sigma)_{ii}, \tau_{ii})
  &=-\frac{\lambda \, h^2}{24 \mu(2\mu+n\lambda)}\sum\limits_{i\not=k=1}^n
   \big(\frac{\partial  \sigma_{kk}}{\partial x_i}, \frac{\partial\tau_{ii}}{\partial x_i}\big)
    +Ch^2\|\sigma\|_2\|\tau_{ii}\|_0.
  \end{split}
  \end{equation}
Since convergence of terms $\frac{\partial\tau_{ii}}{\partial x_i}$ is  unclear,
    we can not obtain directly supercloseness from the previous equation.
The remedy is to use  convergence of divergence
   $\sum\limits_{\ell=1}^n \frac{\partial\tau_{i\ell}}{\partial x_\ell}$ and continuity
 of $\tau_{i\ell}$ across $n-1$ dimensional interior faces
  which are perpendicular to the axis $x_\ell$.
This idea leads to the following decomposition:
\begin{equation}\lab{eq3.19}
  \begin{split}
\sum\limits_{i\not=k=1}^n\bigg(\frac{\partial \sigma_{kk}}{\partial x_i},
    \frac{\partial\tau_{ii}}{\partial x_i}\bigg)
=\sum\limits_{\ell=1}^n\sum\limits_{i\not=k=1}^n
   \bigg(\frac{\partial \sigma_{kk}}{\partial x_i}, \frac{\partial\tau_{i\ell}}{\partial x_\ell}\bigg)
-\sum\limits_{i\not=\ell=1}^n\sum\limits_{i\not=k=1}^n
 \bigg(\frac{\partial \sigma_{kk}}{\partial x_i}, \frac{\partial\tau_{i\ell}}{\partial x_\ell}\bigg).
  \end{split}
  \end{equation}
  Since $\tau_{i\ell}=\sigma_{i\ell,h}-\Pi_{i\ell}\sigma_{i\ell}$,
   the first term on the right--hand side of \eqref{eq3.19} can be estimated
   by the error estimates  presented in \eqref{inter-err5} for the interpolation operator $\Pi_h$
   and convergence  from Theorem \ref{MainError} for the finite element solution $\sigma_h$.
This yields
  \begin{equation}\lab{eq3.20}
  \sum\limits_{\ell=1}^n\sum\limits_{i\not=k=1}^n
  \bigg(\frac{\partial \sigma_{kk}}{\partial x_i}, \frac{\partial\tau_{i\ell}}{\partial x_\ell}\bigg)
  \leq Ch (\|\sigma\|_2^2+\|u\|_1^2).
  \end{equation}
  To analyze the  second term on the right--hand side of \eqref{eq3.19}, we shall explore the continuity to transfer integrations
  on the volume to integrations on the boundary  and  use Lemma \ref{trace}. In fact,
  since the jump $[\tau_{i\ell}]_F$ across face $F$ vanishes for interior face $F\in \mathcal{F}_{x_\ell}$, we have
\begin{equation}\lab{eq3.21}
\begin{split}
\bigg(\frac{\partial \sigma_{kk}}{\partial x_i}, \frac{\partial\tau_{i\ell}}{\partial x_\ell}\bigg)
&=-\bigg(\frac{\partial^2 \sigma_{kk}}{\partial x_i\partial x_\ell}, \tau_{i\ell}\bigg)+\sum\limits_{F\in\mathcal{F}_{x_\ell}}\int_F \big[ \tau_{i\ell} \big]_F \frac{\partial \sigma_{kk}}{\partial x_i} dF\\
&=-\bigg(\frac{\partial^2 \sigma_{kk}}{\partial x_i\partial x_\ell}, \tau_{i\ell}\bigg)+\sum\limits_{F\in\mathcal{F}_{x_\ell}\cap\partial \Omega}\int_F \tau_{i\ell} \frac{\partial \sigma_{kk}}{\partial x_i} dF.
\end{split}
\end{equation}
In order to use  Lemma \ref{trace}, for any $F\in \mathcal{F}_{x_\ell}\cap\partial \Omega$, let $K_F$ be the unique element
 such that $F$ is one of its $n-1$ dimensional faces. Given $v\in L^2(K_F)$, define the constant projection $\Pi_{K_F}^0v$ by
 $$
 \Pi_{K_F}^0v:=\frac{1}{|K_F|}\int_{K_F}vdV.
 $$
This, the trace theorem, inverse estimate and  triangle inequality lead to
\begin{equation}\lab{eq3.21-1}
\begin{split}
\int_F \tau_{i\ell} \frac{\partial \sigma_{kk}}{\partial x_i} dF&=\int_F \tau_{i\ell}(I-\Pi_{K_F}^0)\frac{\partial \sigma_{kk}}{\partial x_i} dF
+ \Pi_{K_F}^0\frac{\partial \sigma_{kk}}{\partial x_i}\int_F \tau_{i\ell}  dF\\
&\leq C\|\tau_{i\ell}\|_{0, K_F}|\sigma_{kk}|_{2, K_F}+Ch^{-1}\|\Pi_{K_F}^0\frac{\partial \sigma_{kk}}{\partial x_i}\|_{0, K_F}\|\tau_{i\ell}\|_{0, K_F}\\
&\leq C \|\tau_{i\ell}\|_{0, K_F}|\sigma_{kk}|_{2, K_F}+Ch^{-1}\|\frac{\partial \sigma_{kk}}{\partial x_i}\|_{0, K_F}\|\tau_{i\ell}\|_{0, K_F}.
\end{split}
\end{equation}
Summing over $F$ in $\mathcal{F}_{x_\ell}\cap\partial \Omega$ and taking in account
the error estimates  presented in \eqref{inter-err5} for the interpolation operator $\Pi_h$ and convergence  from Theorem \ref{MainError} for the finite element solution $\sigma_h$, we arrive at
\begin{equation*}
\begin{split}
\sum\limits_{F\in\mathcal{F}_{x_\ell}\cap\partial \Omega}\int_F \tau_{i\ell} \frac{\partial \sigma_{kk}}{\partial x_i} dF
&\leq C(\|\sigma\|_2+\|u\|_1)(h|\sigma|_{2, \partial\Omega_h}+|\sigma|_{1, \partial\Omega_h}).
\end{split}
\end{equation*}
Hence it follows from Lemma \ref{trace} that
\begin{equation}\lab{eq3.22}
\begin{split}
\sum\limits_{F\in\mathcal{F}_{x_\ell}\cap\partial \Omega}\int_F \tau_{i\ell} \frac{\partial \sigma_{kk}}{\partial x_i} dF
&\leq C(\|\sigma\|_2+\|u\|_1)(h|\sigma|_{2, \Omega}+h^{1/2}\|\sigma\|_{3/2, \Omega}).
\end{split}
\end{equation}
A summary of \eqref{eq3.18} through \eqref{eq3.22}  shows that
\begin{equation}
\sum\limits_{i=1}^n(A(\sigma-\Pi_h\sigma)_{ii}, \tau_{ii})
  \leq Ch^{5/2}(\|\sigma\|_2+\|u\|_1)\|\sigma\|_2.
\end{equation}
Finally, let $\sigma_n=\diag(\sigma_{11}, \cdots, \sigma_{nn})$,
  $\tau_n=\diag(\tau_{11}, \cdots, \tau_{nn})$, $\sigma_s=\sigma-\sigma_n$
  and $\tau_s=\tau-\tau_n$.
 The previous equation, the estimate \eqref{eq3.15} for the shear stress,
estimates \eqref{inter-err4}--\eqref{inter-err5}, and estimates in Theorem \ref{MainError}, yield
\begin{equation*}
(A(\sigma-\Pi_h\sigma), \tau)=(A(\sigma_n-\Pi_h\sigma_n), \tau_n)+(\sigma_s-\Pi_h\sigma_s, \tau_s)
\leq Ch^{5/2}(\|\sigma\|_2+\|u\|_1)\|\sigma\|_2.
\end{equation*}
This completes the proof.
\end{proof}

\begin{theorem} Let $(\sigma, u)$ and $(\sigma_h, u_h)$ be solutions of
    problems \eqref{eqn1} and \eqref{DP},
  respectively.
Suppose that $\sigma \in H^3(\Omega, \mathbb{S})$ and $u\in H^1(\Omega, \mathbb{R}^n)$.
  Then  there holds that
\begin{equation}\lab{eq3.24}
\|\sigma_h-\Pi_h\sigma\|_{H(\d)}^2+\|u_h-P_hu\|_0^2\leq Ch^3(\|\sigma\|_3^2+\|u\|_1^2).
\end{equation}
\end{theorem}
\begin{proof} It is apparent that we can not  derive the desired superconvergence directly
    from  \eqref{eq3.7-1} and \eqref{eq3.16}.
We propose an iterative argument  to show  \eqref{eq3.24}, which consists of the following steps:
\begin{description}

\item[Step 1] 
 By  \eqref{eq3.7-1} and \eqref{eq3.16}, we can deduce the following initial superconvergence result:
\begin{equation}\lab{eq3.25-1}
\|\sigma_h-\Pi_h\sigma\|_0^2+\|u_h-P_hu\|_0^2\leq Ch^{2+\frac{1}{2}}(\|\sigma\|_3^2+\|u\|_1^2).
\end{equation}

\item[Step 2] 
 We  show an intermediate  superconvergence for  $\|\d (\sigma_h-\Pi_h\sigma)\|_0$ 
   based on  \eqref{eq3.1}, which for convenience is recalled as follows
 \begin{equation}\lab{eq3.25}
\begin{split}
&\quad \ \|\Pi_h\sigma-\sigma_h\|_{H(\d)}+\|P_hu-u_h\|_0\\
&\leq C\sup\limits_{0\not=(\tau, v)\in\Sigma_h\times V_h}
  \frac{(A(\sigma-\Pi_h\sigma), \tau)+(u- P_hu, \d \tau)-(\d (\sigma-\Pi_h\sigma), v)
      }{\|\tau\|_{H(\d)}+\|v\|_0}.
\end{split}
\end{equation}
Since the second term on the right--hand side of  \eqref{eq3.25} vanishes and the third term is already analyzed in  \eqref{eq3.5}, we only need
 to show  a better bound for the first term. In fact, it follows
   from  \eqref{eqn1}, \eqref{DP} and \eqref{eq3.25-1} that
 \begin{equation}
 \begin{split}
 (A(\sigma-\Pi_h\sigma), \tau)&=(A(\sigma -\sigma_h), \tau)+(A(\sigma_h-\Pi_h\sigma), \tau)\\
  &=(u_h-u, \d \tau)+(A(\sigma_h-\Pi_h\sigma), \tau)\\
  &=(u_h-P_hu, \d \tau)+(A(\sigma_h-\Pi_h\sigma), \tau)\\
  &\leq Ch^{1+\frac{1}{4}}(\|\sigma\|_3+\|u\|_1)\|\tau\|_{\d}.
\end{split}
 \end{equation}
 Consequently,
 \begin{equation}
 \|\Pi_h\sigma-\sigma_h\|_{H(\d)}\leq Ch^{1+\frac{1}{4}}(\|\sigma\|_3+\|u\|_1).
 \end{equation}
We substitute this estimate into \eqref{eq3.20} to get an improved estimate as
\begin{equation}
  \sum\limits_{\ell=1}^n\sum\limits_{i\not=k=1}^n\bigg(\frac{\partial \sigma_{kk}}{\partial x_i}, \frac{\partial(\sigma_{i\ell,h}-\Pi_{i\ell}\sigma_{i\ell})}{\partial x_\ell}\bigg)
  \leq Ch^{1+\frac{1}{4}} (\|\sigma\|_3^2+\|u\|_1^2).
  \end{equation}
 \item[Step 3] 
  We establish an improved estimate for  the boundary term by 
   putting the estimate $\|\sigma_h-\Pi_h\sigma\|_0$ of \eqref{eq3.25-1} into \eqref{eq3.21-1}:
\begin{equation}
\begin{split}
\sum\limits_{F\in\mathcal{F}_{x_\ell}\cap\partial \Omega}\int_F \tau_{i\ell} \frac{\partial \sigma_{kk}}{\partial x_i} dF
&\leq Ch^{\frac{3}{4}}(\|\sigma\|_3^2+\|u\|_1^2).
\end{split}
\end{equation}
\item[Step 4] 
 We replace those corresponding estimates used in the proof of Lemma \ref{lemma3.8}
    by these improved estimates to obtain
\begin{equation}
(A(\sigma-\Pi_h\sigma), \sigma_h-\Pi_h\sigma_h)\leq Ch^{2+\frac{3}{4}}(\|\sigma\|_3^2+\|u\|_1^2).
\end{equation}
\end{description}
Thus we increase the order of convergence in \meqref{eq3.25-1} from $1+1/4$ to
   $1+1/4+1/8$.
Now we go back to Step 1  and repeat this procedure to get another $1/16$ higher
  order of convergence.
The iteration converges with \meqref{eq3.24}.
\end{proof}

\end{document}